\documentclass{amsart}
\usepackage{graphicx}
\usepackage{tikz}
\usepackage{hyperref}
\usepackage{comment}

\usepackage[T1]{fontenc}

\newtheorem{thm}{Theorem}[section]

\newtheorem{lem}[thm]{Lemma}
\newtheorem{cor}[thm]{Corollary}

\theoremstyle{definition}

\theoremstyle{remark}

\newtheorem{remark}[thm]{Remark}

\numberwithin{equation}{section}

\DeclareMathOperator{\Diff}{Diff}
\DeclareMathOperator{\Homeo}{Homeo}
\DeclareMathOperator{\id}{id}
\DeclareMathOperator{\area}{area}
\DeclareMathOperator{\vol}{vol}

\DeclareMathOperator{\supp}{supp}

\newcommand{\NN}{\mathbf{N}}
\newcommand{\ZZ}{\mathbf{Z}}
\newcommand{\RR}{\mathbf{R}}

\newcommand{\G}{\mathcal{G}}

\newcommand{\x}{\bar{x}}
\newcommand{\z}{\bar{z}}
\newcommand{\REH}{\overline{EH}}
\newcommand{\REG}{\overline{E\Gamma}}

\makeatletter
\@namedef{subjclassname@2020}{\textup{2020} Mathematics Subject Classification}
\makeatother


\begin{document}

\title{Gambaudo--Ghys construction on bounded cohomology}

\author{Mitsuaki Kimura}
\address[Mitsuaki Kimura]{Department of Mathematics, Osaka Dental University, 8-1 Kuzuha Hanazono-cho, Hirakata, Osaka 573-1121, Japan}
\email{kimura-m@cc.osaka-dent.ac.jp}

\begin{abstract}
We consider a generalized Gambaudo--Ghys construction on bounded cohomology and prove its injectivity.
As a corollary, we prove that the third bounded cohomology of the group  of area-preserving diffeomorphisms on the 2-disk is infinite-dimensional. Additionally, we establish similar results for the 2-sphere, the 2-torus, and the annulus.
\end{abstract}

\maketitle

\section{Introduction}

Let $\mathbb{D}$ denote the unit 2-disk equipped with an area form, 
 and $\G$ denote the group $\Diff(\mathbb{D},\partial \mathbb{D},\area)$ of area-preserving diffeomorphisms on $\mathbb{D}$ that are the identity near the boundary $\partial \mathbb{D}$. In \cite{GG}, Gambaudo and Ghys gave a construction of quasimorphisms on $\G$ using the signature of braids. By generalizing their method, Brandenbursky \cite{Bra_knot} defined a linear map $\Gamma \colon Q(P_m) \to Q(\G)$,
where $Q(G)$ denotes the space of homogeneous quasimorphisms on a group $G$, and $P_m$ denotes the pure braid group on $m$ strands.
Let $B_m$ be the braid group on $m$ strands, and $i \colon P_m \to B_m$ the standard inclusion. In \cite{Ish}, Ishida proved that the composition map $\Gamma \circ i^{\ast} \colon Q(B_m) \to Q(\G)$ is injective.
He also proved that the map $\REH_b^2(B_m) \to \REH_b^2(\G)$ induced by $\Gamma \circ i^{\ast}$ is injective, where $\REH_b^n(G)$ denotes the reduced exact bounded cohomology of $G$ (with coefficients in $\RR$).
For definitions on the cohomology and bounded cohomology of groups, see Section \ref{subsec:bddcoh}.
In this paper, we generalize Ishida's result to higher-dimensional bounded cohomology for the case of three strands.
We define a map $\REG_b \colon \REH_b^n(B_m) \to \REH_b^n(\G)$ that generalizes the Gambaudo--Ghys construction and prove the following theorem.
\begin{thm} \label{thm:main}
For $n \geq 2$, the map $\REG_b \colon \REH_b^n(B_3) \to \REH_b^n(\G)$ is injective.
\end{thm}

As a corollary, we obtain the following result.

\begin{cor}\label{cor:disk}
  The dimension of $\REH_b^3(\G)$ is uncountably infinite.
\end{cor}

Our work is inspired by the work of Brandenbursky and Marcinkowski \cite{brandenbursky2019bounded}; 
for a complete Riemannian manifold $M$ of finite volume, under a certain condition on $\pi_1(M)$, they proved that the third bounded cohomology $\REH_b^3(\Diff_0(M,\vol))$ of the
identity component $\Diff_0(M,\vol)$ of the volume-preserving diffeomorphism group on $M$ is uncountably infinite-dimensional.
Note that their result does not yield Corollary \ref{cor:disk} 
since $\pi_1(\mathbb{D})$ is trivial. 
We also note that Nitsche \cite{Nit} has generalized the work of Brandenbursky and Marcinkowski and ours to higher degrees.

We also prove similar results for compact surfaces $\Sigma$ with non-negative Euler characteristic $\chi(\Sigma) \geq 0$. Let $B_m(\Sigma)$ and $P_m(\Sigma)$ denote the braid group and the pure braid group on a surface $\Sigma$, respectively.
Let $\G_{\Sigma}$ denote the identity component of the group $\Diff_{0}(\Sigma,\partial \Sigma ,\area)$ of area-preserving diffeomorphisms on $\Sigma$ that are the identity near the boundary $\partial \Sigma$.
For $m \in \NN$, let $K(\Sigma,m)$ denote the kernel of the forgetful map ${\rm MCG}(\Sigma,m) \to {\rm MCG}(\Sigma)$ (see Section \ref{subsec:braid}).
We also define the map $\REG_b^{\Sigma} \colon \REH_b^n(K(\Sigma,m)) \to \REH_b^n(\G_{\Sigma})$ prove the following.

\begin{thm} \label{thm:S2T2}
Let $\Sigma$ be a compact oriented surface such that $\chi(\Sigma)\geq 0$ and set $m = 2+ \chi(\Sigma)$.
For $n\geq 2$, the map $\REG_b^{\Sigma} \colon \REH_b^n(K(\Sigma,m)) \to \REH_b^n(\G_{\Sigma})$ is injective.
\end{thm}

Note that the case $n=2$ is proved by Brandenbursky and Marcinkowski \cite[Theorem 2.5]{BM_ent}.
We also note that the case $\Sigma=\mathbb{D}$ corresponds to Theorem \ref{thm:main} since $K(\mathbb{D},m)$ is isomorphic to the braid group $B_m$.
Similarly to Corollary \ref{cor:disk}, we obtain the following result.

\begin{cor}\label{cor:S2T2}
  Let $\Sigma$ be a compact oriented surface such that $\chi(\Sigma)\geq 0$.
 The dimension of $\REH_b^3(\G_{\Sigma})$ is uncountably infinite.
\end{cor}

We note that Corollary \ref{cor:S2T2} is not deduced from the result of Brandenbursky and Marcinkowski \cite{brandenbursky2019bounded}.
On the other hand, their result covers the case of surfaces with negative Euler characteristics. Therefore, in some sense, our results and theirs are complementary to each other in the case of 2-manifolds. Namely, we have the following theorem.

\begin{thm}\label{thm:surface}
  For any compact oriented surface $\Sigma$,
 the dimension of $\REH_b^3(\G_{\Sigma})$ is uncountably infinite.
\end{thm}

\section{Preliminary}

\subsection{(Bounded) cohomology of groups} \label{subsec:bddcoh}

We review the definitions on (bounded) cohomology of groups. Let $G$ be a group.
A (homogeneous) \emph{$n$-cochain} on $G$ is a function $c \colon G^{n+1} \to \RR$ such that $c(g_0 h,\dots,g_n h)=c(g_0,\dots,g_n)$ for any $g_0,\dots,g_n ,h \in G$.
Let $C^n(G)$ denote the set of homogeneous $n$-cochains. We define the coboundary map $\delta \colon C^{n-1}(G) \to C^n(G)$ by
\[ \delta c (g_0,\dots , g_n) = \sum_{i=0}^n (-1)^i c(g_0,\dots, \widehat{g_i}, \dots, g_n)  \]
for $c \in C^{n-1}(G)$.
Here, the symbol $\widehat{g}$ means that we omit the entry $g$.
The cochain complex $(C^{\bullet}(G),\delta)$ defines the \emph{group cohomology} $H^{\bullet}(G)$ of $G$.
For a cochain $c \in C^n(G)$, we define $\|c\|_{\infty} \in [0,\infty]$ by
\[ \|c\|_{\infty} = \sup_{g_0,\dots,g_n\in G} |c(g_0,\dots,g_n)|. \]
We say that a cochain $c \in C^n(G)$ is \emph{bounded} if $\|c\|_{\infty} \in [0,\infty)$.
Let $C_b^n(G)$ denote the set of bounded $n$-cochains.
The cochain complex $(C_b^{\bullet}(G),\delta)$ defines the \emph{bounded cohomology} $H_b^{\bullet}(G)$ of $G$.
The inclusion $C_b^n(G) \to C^n(G)$ induces a homomorphism $H_b^n(G) \to H^n(G)$, which is called the \emph{comparison map}. The kernel of the comparison map $H_b^n(G) \to H^n(G)$ is called the \emph{exact bounded cohomology} and is denoted by $EH_b^n(G)$.
The norm $\| \cdot \|_{\infty}$ on $C^n_b(G)$ induces the canonical semi-norm $\|\cdot \|$ on $H_b^n(G)$ defined by  
\[ \| u \| = \inf_{[c]=u} \|c\|_{\infty} \]
for $u \in H_b^n(G)$.
The quotient space of $EH_b^n(G)$ by its norm zero subspace is called the \emph{reduced exact bounded cohomology} and is denoted by $\REH_b^n(G)$.

We summarize several facts which we use later. 

\begin{thm}[{\cite[p.39]{Gromov}}] \label{thm:amenble_galois}
  Let $G$ be a group, $H$ a normal subgroup of $G$, and $i \colon H \to G$ the inclusion map. 
  If $G/H$ is amenable, then the induced map $i^\ast \colon H_b^n(G) \to H_b^n(H)$ is injective and isometric for every $n \geq 1$.
\end{thm}

The following theorem is known as the mapping theorem (for groups).

\begin{thm}[{\cite[p.40]{Gromov}}] \label{thm:mapping}
  If $\phi \colon G_1 \to G_2$ is a surjective group homomorphism with an amenable kernel, then the induced map $\phi^{\ast} \colon H_b^n(G_2) \to H_b^n(G_1)$ is an isometric isomorphism for every $n \geq 1$.
\end{thm}

If $G$ is an amenable group, then its bounded cohomology $H_b^n(G)$  vanishes for every $n \geq 1$. 
On the other hand, the bounded cohomology of non-positive curvature groups tends to be highly non-trivial. For example, the following theorem is known.

\begin{thm}[{\cite[Corollary 6.5]{FPS}}] \label{thm:acyl_hyp}
  If $G$ is an acylindrically hyperbolic group, then the dimension of $\REH_b^3(G)$ is uncountably infinite.
\end{thm}

Examples of acylindrically hyperbolic groups include: 
non-elementary hyperbolic groups, relatively hyperbolic groups, mapping class groups of hyperbolic surfaces, outer automorphism groups of non-abelian free groups, and most 3-manifold groups (see \cite{Osin} for more information on acylindrically hyperbolic groups).


\subsection{Braid groups} \label{subsec:braid}
Let $M$ be a compact connected oriented manifold. Let $X_m(M)$ denote the \emph{\textup{(}orderd\textup{)} configuration space} of $m$ points in $M$, i.e.,
\[ X_n(M) = \{ (x_1,\dots,x_m) \in M^m \mid x_i \neq x_j \; \text{if} \; i \neq j \}. \]
Note that $X_m(M)$ is a codimension 0 submanifold of $M^m$.
The fundamental group of $X_m(M)$ is called the \emph{pure braid group on $m$ strands} on $M$ and denoted by $P_m(M)$.
If $\dim M \geq 3$, it is known that the inclusion $X_m(M) \to M^m$ induces an isomorphism $P_m(M) \to  \pi_1(M^m) \cong \pi_1(M)\times \cdots \times \pi_1(M)$ \cite[Theorem 1.5]{Birman}.
Hence, we are especially interested in the case of $\dim M = 2$.

Let $\Sigma$ be a compact connected oriented 2-dimensional manifold.
The action of the symmetric group 
defines the quotient space $X_m(\Sigma)/\mathfrak{S}_m$, which is called the \emph{unordered configuration space}.
The fundamental group of $X_m(\Sigma)/\mathfrak{S}_m$ is called the \emph{braid group on $m$ strands} on $\Sigma$ and denoted by $B_m(\Sigma)$. There exists a short exact sequence
\[ 1 \to P_m(\Sigma) \to B_m(\Sigma) \to \mathfrak{S}_m \to 1. \]
Thus, we regard $P_m(\Sigma)$ as a normal subgroup of $B_m(\Sigma)$.
Note that $B_m(\mathbb{D})$ is the ordinary Artin braid group $B_m$ and $P_m(\mathbb{D})$ is the pure braid group $P_m$.

Fix a base point $\z$ of the unordered configuration space $X_m(\Sigma)/\mathfrak{S}_m$. Let $P=\{ x_1,\dots, x_m \}$ be a set of distinct $m$ points in $\Sigma$. We define the evaluation map ${\rm ev}_{\z} \colon \Homeo(\Sigma) \to X_m(\Sigma)/\mathfrak{S}_m$ by ${\rm ev}_{\z}(g)=g \cdot \z$, 
where the action of $\Homeo(\Sigma)$ on  $X_m(\Sigma)/\mathfrak{S}_m$ is induced by the diagonal action.
It is known that ${\rm ev}_{\z}$ is a locally trivial fibration with fiber $\Homeo(\Sigma-P)$ (see \cite[Lemma 1.35]{KT}).
Thus, this fibration induces the long exact sequence
\[ \cdots  \pi_1( \Homeo(\Sigma), {\rm id}_\Sigma ) \xrightarrow{{\rm ev}_{\z}^{\ast}} B_m(\Sigma)  \xrightarrow{Push}  {\rm MCG}(\Sigma,m)  \xrightarrow{Forget} {\rm MCG}(\Sigma) \to 1 \]
and the induced map $Forget \colon {\rm MCG}(\Sigma,m) \to {\rm MCG}(\Sigma)$ is called the \emph{forgetful map}, where ${\rm MCG}(\Sigma,m)$ and ${\rm MCG}(\Sigma)$ are the mapping class groups $\pi_0( \Homeo(\Sigma-P) )$ and $\pi_0( \Homeo(\Sigma) )$, respectively.
Let $K(\Sigma,m)$ denote the kernel of the forgetful map. 
The map $Push \colon B_m(\Sigma) \to {\rm MCG}(\Sigma,m) $ is called the \emph{push map}.
Note that $K(\Sigma,m)={\rm Ker}(Forget)={\rm Im}(Push)$.

If $\Sigma$ has non-empty boundary, then $\Homeo(\Sigma)$ is locally contractible \cite{HD}, and thus $K(\Sigma,m)$ is isomorphic to $B_m(\Sigma)$.
For the case where $\Sigma$ is closed, the following result is known.

\begin{thm}[{\cite[Theorem 4.3]{Birman}}] \label{thm:kernel_forget}
   Let $\Sigma$ be a closed oriented surface of genus $g$. If $g\geq 2$, then $K(\Sigma,m)$ is isomorphic to $B_m(\Sigma)$.
  If $g \geq 1, m\geq2$ or $g=0, m\geq 3$, then $K(\Sigma,m)$ is isomorphic to the central quotient $B_m(\Sigma)/Z(B_m(\Sigma))$.
\end{thm}

\section{Generalized Gambaudo--Ghys construction}

In this section, we discuss a generalized Gambaudo--Ghys construction. 

\subsection{The braid $\gamma$}

Set $\G = \Diff(\mathbb{D},\partial\mathbb{D},\area)$ and fix a base point $\bar{z} = (z_1,\dots, z_m) \in X_m(\mathbb{D})$. For simplicity, we assume that $\mathbb{D}$ is equipped with the standard Euclidean area form.
For every $g \in \G$ and almost every $\bar{x} = (x_1,\dots,x_m) \in X_m(\mathbb{D})$, we define a pure braid $\gamma(g,\bar{x}) \in P_m$ as follows.
We take an isotopy $\{ g_t \}_{0 \leq t \leq 1}$ of $g$ such that $g_0 = \id_{\mathbb{D}}$, $g_1 = g$, and $g_t \in \G$ for every $t \in [0,1]$.
We define a loop $l(\{ g_t \},\x) \colon [0,1] \to X_m(\mathbb{D})$ in $X_m(\mathbb{D})$ by
\[ l(\{ g_t \},\x)(t) =
\begin{cases}
\{ (1-3t)z_i + 3tx_i \}_{i=1,\dots,m} & \text{if} \; 0\leq t \leq 1/3,  \\
\{ g_{3t-1}(x_i) \}_{i=1,\dots,m} & \text{if} \; 1/3 \leq t \leq 2/3,  \\
\{ (3-3t)g(x_i) + (3t-2)z_i \}_{i=1,\dots,m} & \text{if} \; 2/3 \leq t \leq 1.
\end{cases}\]
The braid $\gamma(g,\bar{x}) \in P_m$ is defined as the element of $\pi_1(X_m(\mathbb{D}) , \bar{z})$ represented by the loop $l(\{ g_t \},\x)$.
Although $\gamma(g,\bar{x})$ is not defined for every $\x \in X_m(\mathbb{D})$, 
there exists a full measure subspace $\Omega_m$ of $X_m(\mathbb{D})$ with the following property: the braid $\gamma(g,\x)$ is defined if and
only if both $\x$ and $g\cdot \x$ belong to $\Omega_m$ \cite[Section 3.2]{GP}. 

As is well known, $\G$ is contractible (since $\G$ is homotopy equivalent to ${\rm Diff}(\mathbb{D},\partial \mathbb{D})$, which is contractible \cite{Smale}).
Therefore, the above definition of $\gamma(g,\bar{x})$ does not depend on the choice of an isotopy $\{ g_t \}_{0 \leq t \leq 1}$. 

\subsection{The maps $\Gamma_b$ and $\Gamma$} \label{subsec:gamma}

For $c \in C_b^n(B_m)$, we define a map
$\widehat{\Gamma}_b(c) \colon \G^{n+1} \to \RR $ by
\begin{equation} \label{eq:gamma}
\widehat{\Gamma}_b(c)(g_0,\dots,g_n)=
  \int_{\x \in X_m(\mathbb{D})} c( \gamma(g_0,\x),\dots, \gamma(g_n,\x) ) d\x
\end{equation}
for $g_0,\dots,g_n \in \G$. Since $c$ is bounded and the map $\x \mapsto c( \gamma(g_0,\x),\dots, \gamma(g_n,\x) )$ is defined on a full measure subset 
\[\{ \x \in X_m(\mathbb{D}) \mid \x, g_0 \cdot \x, \dots, g_m\cdot \x \in \Omega_m \} = \Omega \cap g_0^{-1}(\Omega_m) \cap \cdots \cap g_n^{-1}(\Omega_m).\]
of $X_m(\mathbb{D})$, the map $\widehat{\Gamma}_b(c)$ is well-defined.

\begin{lem} \label{lem:cochain}
For every $c \in C_b^n(B_m)$, $\widehat{\Gamma}_b(c)$ is a bounded homogeneous cochain.
Moreover, the map $\widehat{\Gamma}_b \colon  C_b^n(B_m) \to  C_b^n(\G)$ is a cochain map.
\end{lem}
\begin{proof}
Since
\[|\widehat{\Gamma}_b(c)(g_0,\dots,g_n)| \leq \vol(X_m(\mathbb{D})) \cdot \| c \|_{\infty}\]
for every $g_0,\dots,g_n \in \G$, $\widehat{\Gamma}_b(c)$ is bounded.
Note that $\gamma(gh , \x)=\gamma(g, h \cdot \x) \gamma(h, \x)$
for $g,h \in \G$, where $\G$ acts diagonally on $X_m(\mathbb{D})$.
Thus,
\begin{align*}
  \widehat{\Gamma}_b(c)(g_0 h ,\dots,g_n h) &=
    \int_{\x \in X_m(\mathbb{D})} c( \gamma(g_0 h, \x),\dots, \gamma(g_n h,\x) ) d\x \\
    &= \int_{\x \in X_m(\mathbb{D})} c( \gamma(g_0,h \cdot  \x)\gamma(h, \x) ,\dots,  \gamma(g_n,h \cdot  \x)\gamma(h, \x) ) d\x \\
    &= \int_{\x \in X_m(\mathbb{D})} c( \gamma(g_0,h \cdot \x) ,\dots, \gamma(g_n,h \cdot \x) ) d\x .
  \end{align*}
Since the action of $h$ preserves the volume form, $\widehat{\Gamma}(c)(g_0 h,\dots,g_n h) = \widehat{\Gamma}(c)(g_0,\dots,g_n)$ and hence $\widehat{\Gamma}(c)$ is homogenous.
By definition, the map $\widehat{\Gamma}$ and the coboundary map $\delta$ commute. Thus, $\widehat{\Gamma}$ is a cochain map.
\end{proof}

By Lemma \ref{lem:cochain}, the map
$\widehat{\Gamma}_b \colon C^n_b(B_m) \to C^n_b(\G)$ induces the homomorphism
$\Gamma_b \colon H^n_b(B_m) \to H^n_b(\G).$

We also define a map $\widehat{\Gamma} \colon C^n(B_m) \to C^n(\G)$ on the ordinary cochain complex as in equation \eqref{eq:gamma}. 
The well-definedness of the map $\widehat{\Gamma}(c) \colon \G^{n+1} \to \RR$ is not obvious since $c \in C^n(B_m)$ is not necessarily bounded, but the map $\widehat{\Gamma}(c)$ is well-defined since the map $\gamma(g, \cdot ) \colon X_m(\mathbb{D}) \to B_m$ has essentially finite image (i.e., there exists a full measure subset of $X_m(\mathbb{D})$ whose image by the map is a finite subset in $B_m$) \cite[Lemma 2.1]{BM_ent}.
Let $\Gamma \colon H^n(B_m) \to H^n(\mathcal{G})$, $E \Gamma_b \colon EH^n_b(B_m) \to EH^n_b(\mathcal{G})$ and $\overline{E\Gamma_b} \colon \overline{EH}^n_b(B_m) \to \overline{EH}^n_b(\mathcal{G})$ be the maps induced by $\widehat{\Gamma}$.

\section{Proof of main result}
In this section, we prove Theorem \ref{thm:main}.
We reduce Theorem \ref{thm:main} to the following key lemma, which corresponds to \cite[Lemma 4.1]{brandenbursky2019bounded}.
Recall that $i \colon P_3 \to B_3$ denotes the inclusion map. 

\begin{lem} \label{lem:key}
There exist a constant $\Lambda > 0$ and a family of homomorphisms $\{\rho_{\epsilon} \colon P_3 \to \G \}_{0<\epsilon <1}$ such that
\[\lim_{\epsilon \to 0} \| \rho_{\epsilon}^{\ast}( \REG_b (u ) ) - \Lambda \cdot i^{\ast}(u)\| = 0 \]
for any $u \in \REH^n_b(B_3)$, where $\rho_{\epsilon}^{\ast} \colon \REH_b^n(\G ) \to \REH_b^n(P_3)$ is the map induced by $\rho_{\epsilon}$.
\end{lem}

Before we prove Lemma \ref{lem:key}, we give the proof of Theorem \ref{thm:main} from Lemma \ref{lem:key}.

\begin{proof}[Proof of Theorem \ref{thm:main}]
Let $u \in \REH_b^n(B_3)$ be a non-trivial class. It means that $\| u \| > 0$, and thus $\|i^{\ast}(u)\|>0$ by Theorem \ref{thm:amenble_galois}. Therefore, by Lemma \ref{lem:key}, we can see that $\|\rho^{\ast}_{\epsilon} (\REG_b (u))\| > 0$ for sufficiently small $\epsilon>0$. It means that $\REG_b (u)$ is non-trivial and hence
$\REG_b $ is injective. 
\end{proof}

Corollary \ref{cor:disk} is deduced from Theorem \ref{thm:main} as follows.

\begin{proof}[Proof of Corollary \ref{cor:disk}]
By Theorem \ref{thm:acyl_hyp}, the dimension of $\REH^3_b(B_3/Z(B_3))$ is uncountably infinite since $B_3/Z(B_3) \cong PSL(2,\ZZ)$ is non-elementary hyperbolic.
The quotient map $B_3 \to B_3/Z(B_3)$ induces an isomorphism $H_b^n(B_3) \to H_b^n(B_3/Z(B_3))$ by Theorem \ref{thm:mapping}. Since $H^3(B_3)=0$ (see \cite[Chapter I]{Va} for example) and $H^3(PSL(2,\ZZ))\cong H^3(\ZZ/2\ZZ) \oplus H^3(\ZZ/3\ZZ) =0$, $EH_b^3(B_3)$ and $EH_b^3(B_3/Z(B_3))$ are also isomorphic.
Therefore, by Theorem \ref{thm:main}, $\REH^3_b(\Diff(\mathbb{D},\area))$ is also uncountably infinite-dimensional.
\end{proof}

In the rest of this section, we prepare for the proof of Lemma \ref{lem:key} in Sections \ref{sec:rho} and \ref{sec:calc_gamma}, and prove Lemma \ref{lem:key} in Section \ref{sec:pf_key_lem}.
The strategy of our proof comes from the work of Brandenbursky and Marcinkowski \cite{brandenbursky2019bounded}, and the method is inspired by the work of Ishida \cite{Ish}.

\subsection{Construction of $\rho_\epsilon$} \label{sec:rho}
\begin{figure}[t]
  \centering
  \includegraphics[width=5cm]{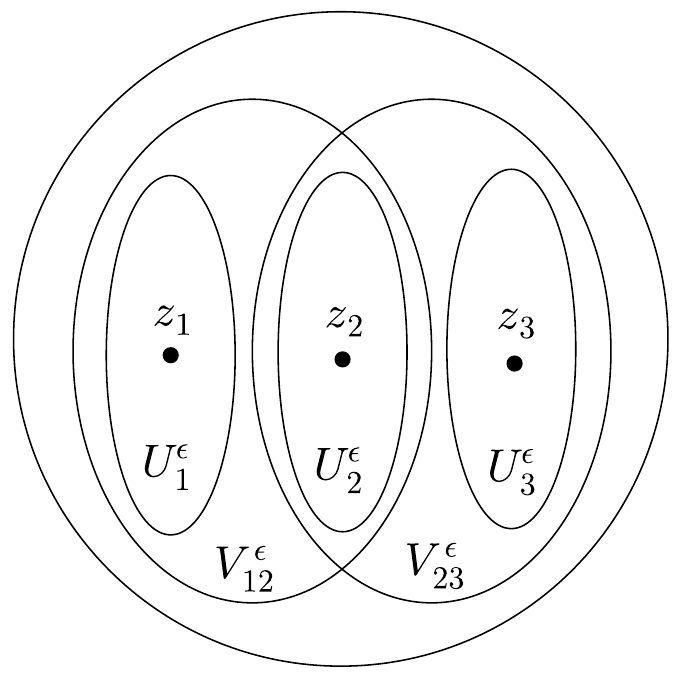}
\caption{Open subsets in $\mathbb{D}$} \label{fig:disk}
\end{figure}

For each $\epsilon$ with $0<\epsilon<1$, 
we construct a homomorphism
$\rho_\epsilon \colon P_3 \to \G$.
Recall that $\z=(z_1,z_2,z_3)$ denotes the base point of $ X_3(\mathbb{D})$.
For each $i=1,2,3$, we take an open neighborhood $U_i^{\epsilon}$ of $z_i$ in $\mathbb{D}$ such that
\begin{itemize}
  \item $U_i^{\epsilon} \cap U_j^{\epsilon} = \emptyset$ if $i\neq j$, and
  \item $\area(U^{\epsilon})=1-\epsilon$, where $U^{\epsilon}=U_1^{\epsilon} \cup U_2^{\epsilon} \cup  U_3^{\epsilon}$.
\end{itemize}
We take open subsets $W^{\epsilon}_{12}$ and $V^{\epsilon}_{12}$ of $\mathbb{D}$ which are diffeomorphic to a disk such that
\begin{itemize}
  \item $U^{\epsilon}_1 \cup U^{\epsilon}_2 \subset \overline{W^{\epsilon}_{12}} \subset V^{\epsilon}_{12}$ and
  \item  $V^{\epsilon}_{12} \cap U^{\epsilon}_3 = \emptyset$.
\end{itemize}
Here, $\overline{W^{\epsilon}_{12}}$ denotes the closure of $W^{\epsilon}_{12}$ in $\mathbb{D}$.
We also take $W^{\epsilon}_{23}$ and $V^{\epsilon}_{23}$ similarly (see Figure \ref{fig:disk}). Finally, we take open disks $W_{123}^{\epsilon}$ and $V_{123}^{\epsilon}$ to be $V_{12}^{\epsilon} \cup V_{23}^{\epsilon} \subset \overline{W_{123}^{\epsilon}} \subset  V_{123}^{\epsilon}$.

We define $\rho_{\epsilon} \colon P_3 \to \G$ as follows.
Set $a_1 = (\sigma_1)^2$, $a_2=(\sigma_2)^2$, and $a_3=\Delta^2$.
Here, $\Delta^2 = (\sigma_1 \sigma_2)^3$ is the full twist.
Then $P_3$ has a presentation
\[P_3 = \langle a_1,a_2,a_3 \mid a_1a_3=a_3a_1, a_2a_3=a_3a_2 \rangle \cong F_2 \times \ZZ.\]
For open disks $V$ and $W$ such that $\overline{W} \subset V$,
let $g_{V,W}\in \G$ denote a diffeomorphism that rotates $W$ once such that $\supp(g_{V,W}) \subset V$.
We define $\rho_{\epsilon} \colon P_3 \to \G$ by
$\rho_{\epsilon}(a_1) = g_{V^{\epsilon}_{12},W^{\epsilon}_{12}}$,
$\rho_{\epsilon}(a_2) = g_{V^{\epsilon}_{23},W^{\epsilon}_{23}}$ and
$\rho_{\epsilon}(a_3) = g_{V^{\epsilon}_{123},W^{\epsilon}_{123}}$.
Note that $\rho_{\epsilon}(a_3)|_{W^{\epsilon}_{123}}=\id_{W^{\epsilon}_{123}}$.
Since $\supp(\rho_{\epsilon}(a_1)) \subset V^{\epsilon}_{12} \subset W^{\epsilon}_{123}$,
$\rho_{\epsilon}(a_1)$ and $\rho_{\epsilon}(a_3)$ commute.
Similarly, $\rho_{\epsilon}(a_2)$ and $\rho_{\epsilon}(a_3)$ are also commutative.
Thus $\rho_{\epsilon}$ is well-defined.

\subsection{Calculation of $\gamma(\rho_{\epsilon}(\alpha),\x)$} \label{sec:calc_gamma}
For $u = [c] \in \REH^n_b(B_3)$, $\rho_{\epsilon}^{\ast}(\REG_b (u)) \in \REH^n_b(P_3)$ is the cohomology class of the cochain defined by
\[(\alpha_0,\dots,\alpha_n) \mapsto \int_{\x \in X_3(\mathbb{D})} c(\gamma(\rho_{\epsilon}(\alpha_0),\x),\dots,  \gamma(\rho_{\epsilon}(\alpha_n),\x) )d\x \]
for $\alpha_0,\dots, \alpha_n \in P_3$.
We calculate $\gamma(\rho_{\epsilon}(\alpha),\x) \in P_3 $
for $\alpha \in P_3$ and $\x=(x_1,x_2,x_3) \in X_3(\mathbb{D})$. To describe it, we prepare several notions.
We say that $x \in X_3(\mathbb{D})$ is an \emph{$\epsilon$-good point} if all of $x_1$, $x_2$, and $x_3$ are in $U^{\epsilon}$. Otherwise, we say that $\x$ is an \emph{$\epsilon$-bad point}.
We say that an $\epsilon$-good point $\bar{x}$ is \emph{of type $(p,q,r)$} if
\[ \# (U^\epsilon_1 \cap \{ x_1, x_2, x_3\}) = p, \quad \# (U^\epsilon_2 \cap \{ x_1, x_2, x_3\}) = q, \quad \# (U^\epsilon_3 \cap \{ x_1, x_2, x_3\}) = r. \]

We define homomorphisms $s_i \colon P_3 \to \ZZ$ ($i=1,2,3$) by $s_i( a_j )=\delta_{ij}$ for $1\leq i,j \leq 3$, where $\delta_{ij}$ is the Kronecker delta. 
For each type $(p,q,r)$, we define a homomorphism $\phi_{pqr} \colon P_3 \to P_3$ by
\begin{equation} \label{eq:braid_disk}
  \phi_{pqr}(\alpha) =
  \begin{cases}
     \alpha & \text{type $(1,1,1)$}, \\
     a_3^{s_1(\alpha)+s_3(\alpha)} & \text{type $(3,0,0)$ or $(2,1,0)$}, \\
     a_3^{s_2(\alpha)+s_3(\alpha)} & \text{type $(0,0,3)$ or $(0,1,2)$}, \\
     a_3^{s_1(\alpha)+s_2(\alpha)+s_3(\alpha)} & \text{type $(0,3,0)$}, \\
     a_1^{s_1(\alpha)}a_3^{s_3(\alpha)} & \text{type $(2,0,1)$}, \\
     a_2^{s_2(\alpha)}a_3^{s_3(\alpha)} & \text{type $(1,0,2)$}, \\
     a_1^{s_1(\alpha)}a_3^{s_2(\alpha)+s_3(\alpha)} & \text{type $ (0,2,1)$}, \\
     a_2^{s_2(\alpha)}a_3^{s_1(\alpha)+s_3(\alpha)} & \text{type $ (1,2,0)$}.
  \end{cases}
\end{equation}
The following is the key to the proof of Lemma \ref{lem:key} (compare with Ishida \cite[Theorem 1.2]{Ish}).

\begin{lem}\label{lem:main_obs}
For almost every $\epsilon$-good point $\x \in X_m(\mathbb{D})$ of type $(p,q,r)$, there exists a braid $\beta(\x) \in B_3$ such that
\[\gamma( \rho_{\epsilon}(\alpha, \x) ) = \beta(\x)\phi_{pqr}(\alpha)\beta(\x)^{-1} \]
for every $\alpha \in P_3$.
\end{lem}

\begin{figure}[t]
  \begin{minipage}[b]{0.45\linewidth}
      \centering
     \includegraphics[width=5.5cm]{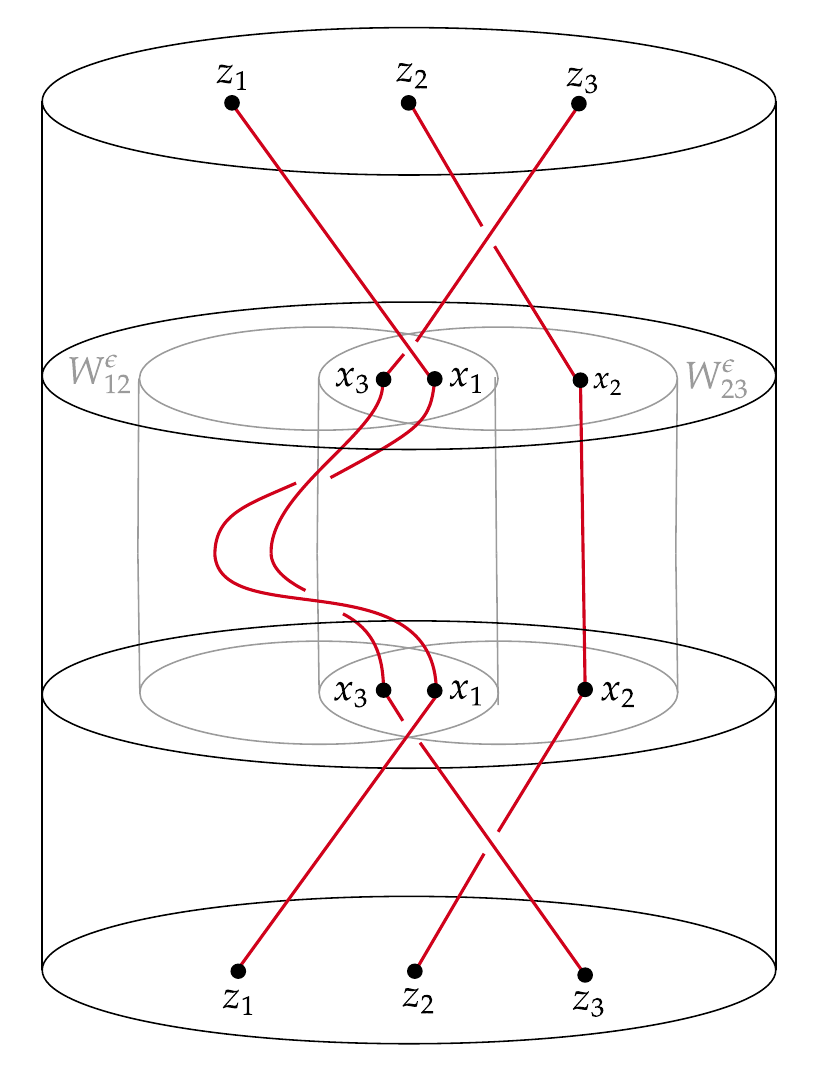}
   \end{minipage}
  \begin{minipage}[b]{0.45\linewidth}
    \centering
     \includegraphics[width=5.5cm]{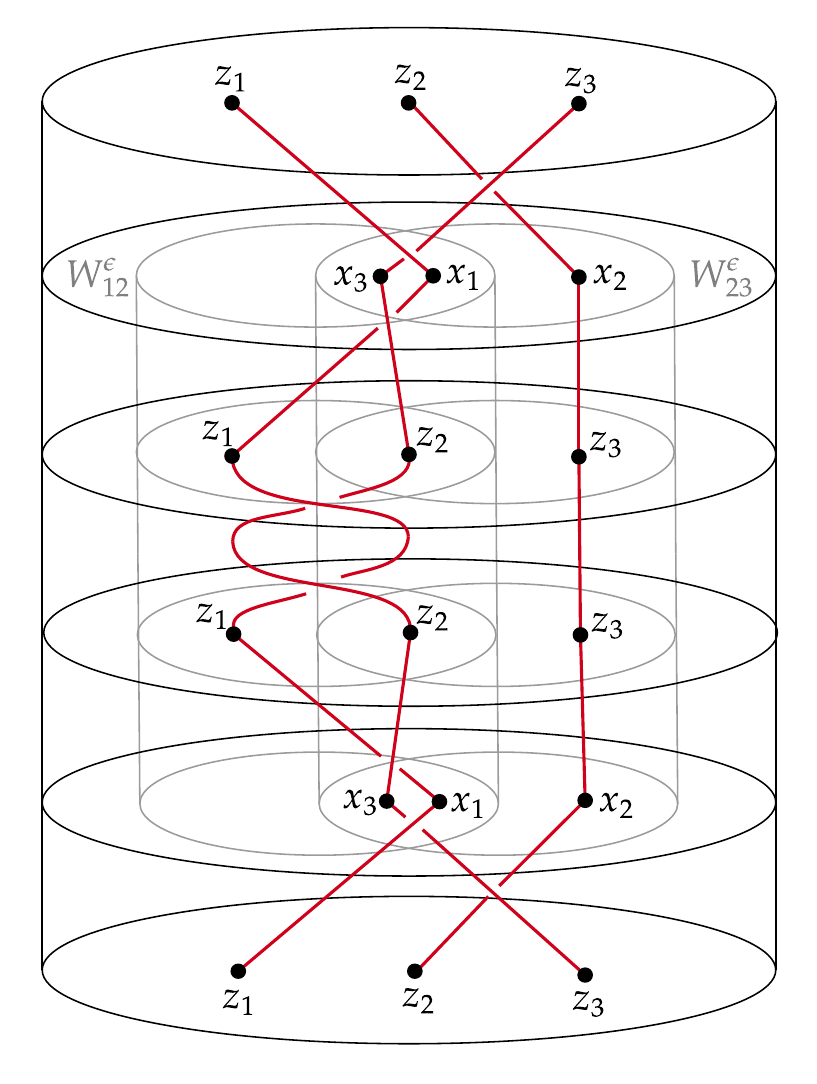}
   \end{minipage}
   \caption{Braids $\gamma(\rho_{\epsilon}(a_1),\x)$ and $\beta(\x) a_1 \beta(\x)^{-1}$ for $\x$ is of type (0,1,2)} \label{fig:calc}
\end{figure}

\begin{proof}
Let $\x=(x_1,x_2,x_3)$ be of type $(p,q,r)$. 
Assume that $x_{i_1},\dots,x_{i_p} \in U^{\epsilon}_1$, 
$x_{j_1},\dots,x_{j_q} \in U^{\epsilon}_2$, and 
$x_{k_1},\dots,x_{k_r} \in U^{\epsilon}_3$. 
The subscript $i$ is defined so that $i_1 < \cdots < i_p$; the same applies to $j$ and $k$. 
We define $\sigma_{\x} \in \mathfrak{S}_3$ by 
\[ \sigma_{\x}=
  \left(
  \begin{matrix}
    i_1 & \cdots & i_p & j_1 & \cdots & j_q & k_1 & \cdots & k_r \\
    1 & \cdots & p & p+1 & \cdots & p+q & p+q+1 & \cdots & p+q+r \\
  \end{matrix}
  \right).
\]
For example, if $x_1,x_3 \in U_2^{\epsilon}$ and $x_2 \in U_3^{\epsilon}$, then $j_1=1$, $j_2=3$, $k_1=2$ and thus
$\sigma_{\x} = \left(
  \begin{matrix}
    1 & 3 & 2 \\
    1 & 2 & 3 \\
  \end{matrix}
  \right).$ Note that $\sigma_{\x}=e$ if $\x$ is of type $(0,0,3)$, $(0,3,0)$ or $(3,0,0)$.

We define $\beta(\bar{x}) \in B_3$ as the element of $\pi_1(X_3(\mathbb{D})/\mathfrak{S}_3 , \bar{z})$ represented by the loop $l \colon [0,1] \to X_3(\mathbb{D})/\mathfrak{S}_3$ defined by
\[ l(t) =
\begin{cases}
\{ (1-2t)z_i + 2tx_i \}_{i=1,2,3} & \text{if} \quad 0\leq t \leq 1/2,  \\
\{ (2-2t)x_i + (2t-1)z_{\sigma_{\x}(i)} \}_{i=1,2,3} & \text{if} \quad 1/2 \leq t \leq 1.
\end{cases}\]
The braid $\beta(\x)$ is defined for almost every $\epsilon$-good point $\x$. Note that the projection $B_3 \to \mathfrak{S}_3$ maps $\beta(\x)$ to $\sigma_{\x}$.
Then, the calculation of $\gamma(\rho_{\epsilon}(a_i),\x)$ for the generators $a_1$, $a_2$, $a_3$ of $P_3$ is as follows (see also Figure \ref{fig:calc}):
\begin{itemize}
    \item $\gamma(\rho_{\epsilon}(a_1),\x) = 
\begin{cases} 
e & \text{if} \quad p+q \leq 1, \\
\beta(\x) a_1 \beta(\x)^{-1} & \text{if} \quad p+q = 2, \\
\beta(\x) a_3 \beta(\x)^{-1}  & \text{if} \quad p+q = 3. 
\end{cases}$

    \item $\gamma(\rho_{\epsilon}(a_2),\x) = 
\begin{cases} 
e & \text{if} \quad q+r \leq 1, \\
\beta(\x) a_2 \beta(\x)^{-1} & \text{if} \quad q+r = 2, \\
\beta(\x) a_3 \beta(\x)^{-1} & \text{if} \quad q+r = 3. 
\end{cases}$
\item $\gamma(\rho_{\epsilon}(a_3),\x)= \beta(\x) a_3 \beta(\x)^{-1}$
\end{itemize}

If $(p,q,r)=(1,1,1)$, since $\gamma(\rho_{\epsilon}(a_i),\x)=\beta(\x)a_i\beta(\x)^{-1}$ for $i=1,2,3$, it follows that $\gamma(\rho_{\epsilon}(\alpha),\x)=\beta(\x)\alpha\beta(\x)^{-1}$ for every $\alpha \in P_3$.
If $(p,q,r)\neq(1,1,1)$, noting that $a_3$ commutes with any braid, we obtain the assertion.
\end{proof}

\subsection{Proof of the key lemma} \label{sec:pf_key_lem}

Now we are ready to prove Lemma \ref{lem:key}.

\begin{proof}[Proof of Lemma \ref{lem:key}]

For each $\epsilon$ with $0<\epsilon<1$, we take an open neighborhood $U_i^{\epsilon}$ of $z_i$ in $\mathbb{D}$ for $i=1,2,3$, and construct the homomorphism $\rho_\epsilon \colon P_3 \to \G$ as in Section \ref{sec:rho}. 
Let $X^{\epsilon}_{pqr}$ denote the set of $\epsilon$-good points in $X_3(\mathbb{D})$ of type $(p,q,r)$ and $Y^{\epsilon}$ denote the set of $\epsilon$-bad points.
We define cochains $c_{pqr}^{\epsilon}, c_Y^{\epsilon} \in C^n_b(P_3)$ by
\begin{align*}
  c_{pqr}^{\epsilon}(\alpha_0,\dots,\alpha_n) & = \int_{\x \in X^{\epsilon}_{pqr}}
  c( \gamma(\rho_{\epsilon}(\alpha_0),\x),\dots,  \gamma(\rho_{\epsilon}(\alpha_n),\x) )d\x,  \\
  c_Y^{\epsilon}(\alpha_0,\dots,\alpha_n) & = \int_{\x \in Y^{\epsilon}}
  c( \gamma(\rho_{\epsilon}(\alpha_0),\x),\dots,  \gamma(\rho_{\epsilon}(\alpha_n),\x) )d\x
\end{align*}
for $\alpha_0,\dots,\alpha_n \in P_3$. Note that
\[\rho_{\epsilon}^{\ast}(\REG_b (u) ) =
 \sum_{p,q,r}[c_{pqr}^{\epsilon}] + [c_Y^{\epsilon}] \in \REH_b^n(P_3).\]

For $c\in C^n(B_3)$ and $\beta \in B_3$, let $\beta \cdot c \in C^n(B_3)$ denote the cochain defined by
\[(\beta \cdot c) (\gamma_0 ,\dots, \gamma_n) = c (\beta \gamma_0 \beta^{-1},\dots, \beta \gamma_n \beta^{-1}).\]
for $\gamma_0,\dots,\gamma_n \in B_3$. 
By Lemma \ref{lem:main_obs}, $c^\epsilon_{pqr}$ satisfies that
\begin{align*}
  c_{pqr}^{\epsilon}(\alpha_0,\dots,\alpha_n) & = \int_{\x \in X^{\epsilon}_{pqr}}
  c( \beta(\x) \phi_{pqr}(\alpha_0) \beta(\x)^{-1},\dots, \beta(\x)\phi_{pqr}(\alpha_n)\beta(\x)^{-1})d\x \\
  &= \sum_{\beta \in B_3} \vol \left( \{ \x \in X_{pqr}^{\epsilon} \mid \beta(\x)=\beta \} \right) (\beta \cdot c)(\phi_{pqr}(\alpha_0) ,\dots, \phi_{pqr}(\alpha_n))
\end{align*}
for any $\alpha_0,\dots,\alpha_n \in P_3$. Since $[\beta \cdot c]=[c]=u$ for any $\beta \in B_3$, it holds that
\begin{equation} \label{eq:c_pqr}
  [c_{pqr}^{\epsilon}] = \vol (X^{\epsilon}_{pqr}) \cdot \phi^{\ast}_{pqr}(i^{\ast}(u)).
\end{equation}

If $(p,q,r)= (1,1,1)$, since $\phi_{111}=\id$ and by equation \eqref{eq:c_pqr},
\[[c_{111}^{\epsilon}] = \vol (X^{\epsilon}_{111}) \cdot i^{\ast}(u) =
3! \cdot \area(U^{\epsilon}_1)\area(U^{\epsilon}_2)\area(U^{\epsilon}_3) \cdot  i^{\ast} (u).\]

If $(p,q,r)\neq (1,1,1)$, by equation \eqref{eq:braid_disk}, the homomorphism $\phi_{pqr}$ factors through
the abelian subgroup $\langle a_i, a_3 \rangle \cong \ZZ^2$ of $P_3$, where $i=1$ or $i=2$.
Since $\ZZ^2$ is amenable, $\REH_b^n(\ZZ^2)=0$.
Thus $\phi_{pqr}^{\ast} = 0 $ and hence $[c^{\epsilon}_{pqr}]=0$ by equation \eqref{eq:c_pqr}.

By the definition of $c_Y^{\epsilon}$,
\[ |c_Y^{\epsilon}(\alpha_0,\dots,\alpha_n)| \leq \vol(Y^{\epsilon}) \| c \|_{\infty}. \]
Since $\vol(Y^{\epsilon}) = \vol(X_3(\mathbb{D})) -\vol(U^{\epsilon}\times U^{\epsilon} \times U^{\epsilon}) = 1 -(1-\epsilon)^3 $,
$\lim_{\epsilon \to 0} \|[c_Y^{\epsilon}]\| = 0$.

Therefore, by setting $\Lambda = \lim_{\epsilon \to 0} 3! \cdot  \area(U^{\epsilon}_1)\area(U^{\epsilon}_2)\area(U^{\epsilon}_3)$,
\[\lim_{\epsilon \to 0} \| \rho_{\epsilon}^{\ast}(\REG_b (u)) - \Lambda \cdot  i^{\ast} (u) \| = 0. \qedhere\]
\end{proof}

\section{The case of other surfaces} \label{sec:other_surfaces}
In this section, 
 we apply the argument from the previous section
to the case of other surfaces and prove Theorem \ref{thm:S2T2}.
Let $\Sigma$ be a compact surface with an area form. We set $\G_{\Sigma} = \Diff_{0}(\Sigma,\partial\Sigma,\area)$.

We provide a generalized Gambaudo--Ghys construction on surfaces (see also \cite[Section 2]{BM_ent}). 
Take a continuous map $\iota \colon \mathbb{D} \to \Sigma$ such that $\iota|_{\mathbb{D} \setminus \partial\mathbb{D}}$ is injective and $\iota(\mathbb{D} \setminus \partial\mathbb{D})$ is of full measure in $\Sigma$.
Take a base point $\bar{z} = (z_1,\dots,z_m)$ of $X_m(\Sigma)$ so that $z_i \in \iota(\mathbb{D} \setminus \partial\mathbb{D})$ for each $i$.
Let $g \in \G_{\Sigma}$ and fix an isotopy $\{g_t\}_{0\leq t \leq 1}$ of $g$.
For $\x =(x_1,\dots,x_m) \in \iota_{\ast}(\Omega_m)$, we define the loop $l(\{g_t\},\x) \colon [0,1] \to X_m(\Sigma)$ by 

\[ l(\{ g_t \},\x)(t) =
\begin{cases}
\left\{ \iota\Bigl
( (1-3t) \cdot  \iota^{-1}(z_i) + 3t \cdot \iota^{-1}(x_i) \Bigr) \right\}_{i=1,\dots,m} & \text{if} \; 0\leq t \leq 1/3,  \\
\{ g_{3t-1}(x_i) \}_{i=1,\dots,m} & \text{if} \; 1/3 \leq t \leq 2/3,  \\
\left\{ \iota\Bigl( (3-3t)\cdot \iota^{-1}( g(x_i)) + (3t-2) \cdot \iota^{-1}(z_i) \Bigr) \right\}_{i=1,\dots,m} & \text{if} \; 2/3 \leq t \leq 1.
\end{cases}\]
Let $\gamma(\{g_t\},\x)$ denote an element of $\pi_1(X_m(\Sigma),\bar{z}) \cong P_m(\Sigma)$ represented by the loop $l(\{g_t\},\x)$.
In general, $\gamma(\{g_t\},\x)$ depends on the choice of an isotopy $\{g_t\}$ but $Push(\gamma(\{g_t\},\x)) \in \mathrm{MCG}(\Sigma,m)$ does not: 
if $\{g'_t\}_{0\leq t \leq 1}$ is another isotopy of $g$, then $\gamma(\{g_t\},\x)\gamma(\{g'_t\},\x)^{-1} \in \mathrm{Im}(ev_{z}^{\ast})=\mathrm{Ker}(Push)$. 
Thus, $\gamma(\{g_t\},\x)$ defines an element of $K(\Sigma,m)$ and we write this element as $\gamma(g,\x)$.

In this way, we can define the map
$\widehat{\Gamma}^{\Sigma}_b \colon C_b^n(K(\Sigma,m)) \to C_b^n(\G_{\Sigma})$
in the same way as in Section \ref{subsec:gamma} since $\iota_{\ast}(\Omega_m)$ is a full measure subset of $X_m(\Sigma)$.
Here, $\iota_\ast \colon X_m(\mathbb{D}) \to X_m(\Sigma)$ is the map induced by $\iota$.
This map $\widehat{\Gamma}_b^{\Sigma}$ is a cochain map by the same arguments as in Lemma \ref{lem:cochain}, and induces the map $\Gamma_b^{\Sigma} \colon H_b^n(K(\Sigma,m)) \to H_b^n(\G_{\Sigma})$.
Moreover,
the map
$\widehat{\Gamma}^{\Sigma} \colon C^n(K(\Sigma,m)) \to C^n(\G_{\Sigma})$ defined by
\[ \widehat{\Gamma}^{\Sigma}(c)(g_0,\dots,g_n)=
  \int_{\x \in X_m(\Sigma)} c( \gamma(g_0,\x),\dots, \gamma(g_n,\x) ) d\x\]
  is well-defined since the map $\gamma(g, \cdot ) \colon X_m(\Sigma) \to MCG(\Sigma,m)$ has essentially finite image \cite[Lemma 2.1]{BM_ent}.
Then $\widehat{\Gamma}^{\Sigma}$ induces the map $\Gamma^{\Sigma} \colon H^n(K(\Sigma,m)) \to H^n(\G_{\Sigma})$ and hence induces the map $\REG_b^{\Sigma} \colon \REH_b^n(K(\Sigma,m)) \to \REH_b^n(\G_{\Sigma})$.

\subsection{The case of an annulus}

Let $\mathbb{A}$ denote an annulus $S^1 \times [0,1]$.
The braid group $B_m(\mathbb{A})$ on $\mathbb{A}$ is isomorphic to the inverse image $\pi^{-1}(\mathfrak{S}_{m})$ of the subgroup $\mathfrak{S}_{m} \subset \mathfrak{S}_{m+1}$ of $\mathfrak{S}_{m+1}$ by the projection $\pi \colon B_{m+1} \to \mathfrak{S}_{m+1}$ \cite[Theorem 2]{MR1902362}. 
This is because a ``pillar'' in $\mathbb{A} \times [0,1]$ can be seen as a ``fixed'' strand (Figure \ref{fig:annular_braid}, see also \cite[Section 2]{MR1902362}). 
Namely, the pure braid group $P_m(\mathbb{A})$ on $\mathbb{A}$ is isomorphic to the ordinary pure braid group $P_{m+1}$, thus we identity them. 
Let $i \colon P_m(\mathbb{A}) \to B_m(\mathbb{A})$ be the inclusion map and set $j = Push \circ i \colon P_m(\mathbb{A}) \to K(\mathbb{A},m)$.

\begin{lem} \label{lem:annulus}
There exist a constant $\Lambda > 0$ and a family of homomorphisms $\{\rho_{\epsilon} \colon P_2(\mathbb{A}) \to \G_{\mathbb{A}} \}_{0<\epsilon <1}$ such that
\[\lim_{\epsilon \to 0} \| \rho_{\epsilon}^{\ast}( \REG^{\mathbb{A}}_b (u ) ) - \Lambda \cdot j^{\ast}(u)\| = 0 \]
for any $u \in \REH^n_b(K(\mathbb{A},2))$.
\end{lem}

\begin{proof}
  For each $\epsilon$, we take open neighborhood $U_i^{\epsilon}$ of $z_i \in \mathbb{A}$ $(i=1,2)$ so that
  \begin{itemize}
    \item $U_1^{\epsilon} \cap U_2^{\epsilon} = \emptyset$ and
    \item $\area(U^{\epsilon})=1-\epsilon$, where $U^{\epsilon}=U^{\epsilon}_1 \cup U^{\epsilon}_2$.
  \end{itemize}
  Moreover, we take open subsets $W^{\epsilon}_1$ and $V^{\epsilon}_1$ of $\mathbb{A}$ that are diffeomorphic to an annulus so that
  \begin{itemize}
    \item $U_1^{\epsilon} \subset \overline{W^{\epsilon}_{1}} \subset V^{\epsilon}_{1}$,
    \item $V^{\epsilon}_1 \cap U_2^{\epsilon}  = \emptyset$ and
    \item the inclusion map $W^{\epsilon}_{1} \to \mathbb{A}$ induces an isomorphism $\pi_1(W^{\epsilon}_{1}) \to \pi_1(\mathbb{A})$.
  \end{itemize}
  We also take $W^{\epsilon}_{2}$ and $V^{\epsilon}_{2}$ in a similar way (Figure \ref{fig:annulus}).
  Finally, we take open annulus $W_{12}^{\epsilon}$ and $V_{12}^{\epsilon}$ to be $V_{1}^{\epsilon} \cup V_{2}^{\epsilon} \subset \overline{W_{12}^{\epsilon}} \subset  V_{12}^{\epsilon}$.

  \begin{figure}[t]
  \begin{tabular}{cc}
   \begin{minipage}{0.45\linewidth}
     \centering
     \includegraphics[width=5cm]{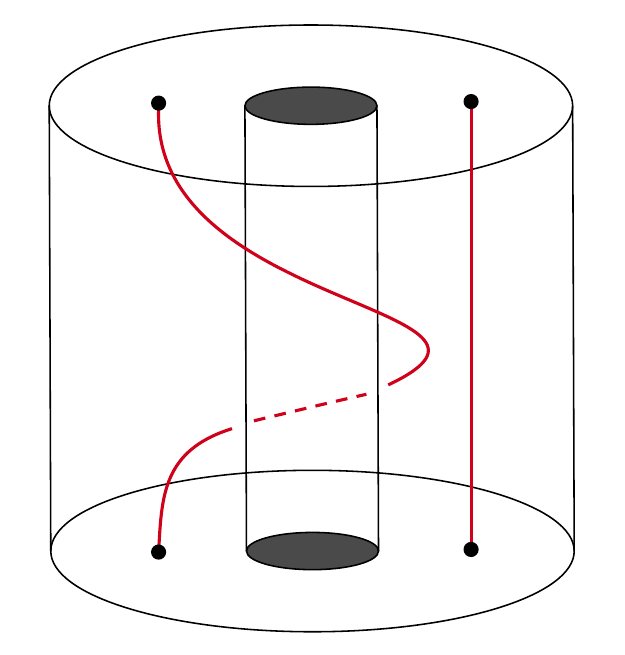}
   \caption{The $2$-braid $(\sigma_1)^2$ on $\mathbb{A}$}
   \label{fig:annular_braid}
   \end{minipage}
   \begin{minipage}{0.45\linewidth}
     \centering
     \includegraphics[width=5cm]{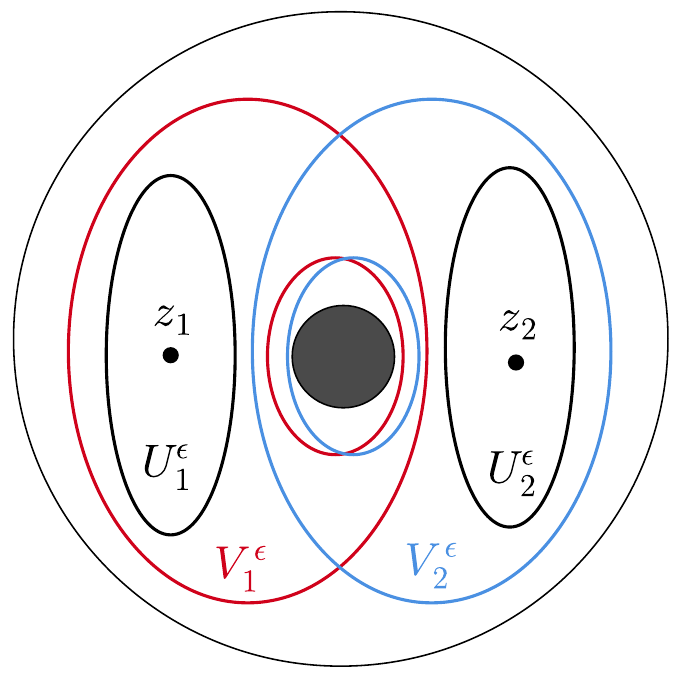}
   \caption{Open subsets in $\mathbb{A}$}
   \label{fig:annulus}
   \end{minipage}
    \end{tabular}
  \end{figure}

  We define $\rho_{\epsilon} \colon P_2(\mathbb{A}) \to \G_{\mathbb{A}}$ as follows.
  Recall that $P_2(\mathbb{A}) \cong P_3$ has a presentation
  \[P_3 = \langle a_1,a_2,a_3 \mid a_1a_3=a_3a_1, a_2a_2=a_3a_2 \rangle \cong F_2 \times \ZZ, \]
  where $a_1=(\sigma_1)^2$, $a_2=(\sigma_2)^2$, and $a_3=\Delta^2=(\sigma_1\sigma_2)^3$.
  For open annuli $V$ and $W$ such that $\overline{W} \subset V$,
  let $g_{V,W}\in \G_{\mathbb{A}}$ denote a diffeomorphism which rotates $W$ once such that $\supp(g_{V,W}) \subset V$.
  We define $\rho_{\epsilon} \colon P_2(\mathbb{A}) \to \G_{\mathbb{A}}$ by
  $\rho_{\epsilon}(a_1) = g_{V^{\epsilon}_{1},W^{\epsilon}_{1}}$,
  $\rho_{\epsilon}(a_2) = g_{V^{\epsilon}_{2},W^{\epsilon}_{2}}$ and
  $\rho_{\epsilon}(a_3) = g_{V^{\epsilon}_{12},W^{\epsilon}_{12}}$.


We say that $\x =(x_1,x_2) \in X_2(\mathbb{A})$ is an \emph{$\epsilon$-good point} if both $x_1$ and $x_2$ are in $U^{\epsilon}$.
We say that an $\epsilon$-good point $\bar{x}$ is \emph{of type $(p,q)$} if
\[ \# (U^\epsilon_1 \cap \{ x_1, x_2\}) = p, \quad \# (U^\epsilon_2 \cap \{ x_1, x_2\}) = q. \]

  Let $\x\in X_2(\mathbb{A})$ be an $\epsilon$-good point of type $(p,q)$.
  If $(p,q) \neq (1,1)$, we can see that $\gamma(\rho_{\epsilon}(\alpha), \x ) \in Z(P_2(\mathbb{A})) = \langle \Delta^2 \rangle $ for any $\alpha \in P_2(\mathbb{A})$.
  By an argument similar to the proof of Lemma \ref{lem:key}, we can prove that
  \[\lim_{\epsilon \to 0} \| \rho_{\epsilon}^{\ast}(\REG_b^{\mathbb{A}}(u)) - \Lambda \cdot  j^{\ast} (u) \| = 0\]
  by setting
  $\Lambda = \lim_{\epsilon \to 0} 2! \cdot \area(U^{\epsilon}_1) \area(U^{\epsilon}_2)$.
\end{proof}

\subsection{The case of a sphere} \label{sec:sphere}
Let $\mathbb{S}$ denote the 2-sphere.
We summarize several facts on the braid group on $\mathbb{S}$ 
(see \cite[Section 4.1]{GJP} for example).

An inclusion $\mathbb{D} \to \mathbb{S}$ induces the projection $B_m \to B_m(\mathbb{S})$, and let $\delta_i$ denote the image of $\sigma_i$ by this projection. It is known that the kernel of this projection is normally generated by $\sigma_1\sigma_2\cdots\sigma_{m-2}\sigma_{m-1}^2\sigma_{m-2}\cdots\sigma_2\sigma_1$.
The natural map $X_{m-1}(\mathbb{D}) \to X_{m}(\mathbb{S})$ induces
the map $P_{m-1} \to P_m(\mathbb{S})$ and it is known to be surjective.
If $m \geq 4$, $Z(P_m(\mathbb{S}))=Z(B_m(\mathbb{S}))$ is generated by the full twist $\xi^2 = (\delta_1\delta_2 \cdots \delta_{m-1})^m$ and $\xi^2$ has order two.

We consider in particular the case $m=4$.
Then $P_4(\mathbb{S})$ has a presentation
\[ P_4(\mathbb{S}) = \langle a_1,a_2,a_3 \mid a_1a_3=a_3a_1, a_2a_3=a_3a_2, (a_3)^2=e \rangle \cong F_2 \times \ZZ/2\ZZ, \]
where $a_1= (\xi_1)^2$, $a_2= (\xi_2)^2$, and $a_3= {\xi}^2$.
Let $i \colon P_4(\mathbb{S}) \to B_4(\mathbb{S})$ be the inclusion.
By Theorem \ref{thm:kernel_forget} and since $Z(P_4(\mathbb{S}))=Z(B_4(\mathbb{S}))$,
the map $Push \circ i \colon P_4(\mathbb{S}) \to K(\mathbb{S},4)$ induces
the map $j \colon P_4(\mathbb{S})/ Z(P_4(\mathbb{S})) \to  K(\mathbb{S},4)$.
For an element $\alpha \in G$ of a group $G$, let $\overline{\alpha} \in
G/Z(G)$ denote the equivalence class of $\alpha$.
We regard the group $P_4(\mathbb{S})/Z(P_4(\mathbb{S}))$ as a free group $F_2$ of rank 2 generated by $\overline{a_1}$ and $\overline{a_2}$.

\begin{lem} \label{lem:sphere}
There exist a constant $\Lambda > 0$ and a family of homomorphisms $\{\rho_{\epsilon} \colon F_2 \to \G_{\mathbb{S}} \}_{0<\epsilon <1}$  such that
\[\lim_{\epsilon \to 0} \| \rho_{\epsilon}^{\ast}( \REG_b^{\mathbb{S}} (u ) ) - \Lambda \cdot j^{\ast}(u)\| = 0 \]
for any $u \in \REH^n_b(K(\mathbb{S},4))$.
\end{lem}

\begin{proof}
  For each $\epsilon$, we take open neighborhoods $U_i^{\epsilon}$ of $z_i \in \mathbb{S}$ $(i=1,2,3,4)$ so that
  \begin{itemize}
    \item $U_i^{\epsilon} \cap U_j^{\epsilon} = \emptyset$ if $i\neq j$ and
    \item $\area(U^{\epsilon})=1-\epsilon$, where $U^{\epsilon}=U_1^{\epsilon} \cup U_2^{\epsilon} \cup  U_3^{\epsilon} \cup U_4^{\epsilon}$.
  \end{itemize}
  Moreover, we take open subsets $W^{\epsilon}_{12}$ and $V^{\epsilon}_{12}$ of $\mathbb{S}$ which are diffeomorphic to a disk so that
  \begin{itemize}
    \item $U^{\epsilon}_1 \cup U^{\epsilon}_2 \subset \overline{W^{\epsilon}_{12}} \subset V^{\epsilon}_{12}$,
    \item  $V^{\epsilon}_{12} \cap U^{\epsilon}_3 = \emptyset$ and
    \item  $V^{\epsilon}_{12} \cap U^{\epsilon}_4 = \emptyset$.
  \end{itemize}
  We also take $W^{\epsilon}_{23}$ and $V^{\epsilon}_{23}$ similarly (see Figure \ref{fig:sphere}).
  We define $\rho_{\epsilon} \colon F_2 \to \G_{\mathbb{S}}$ as in the case of the disk. 
  We define $s_i \colon F_2 \to \ZZ$ by $s_i(\overline{a_j})=\delta_{ij}$.
  \begin{figure}[t]
    \centering
    \includegraphics[width=5cm]{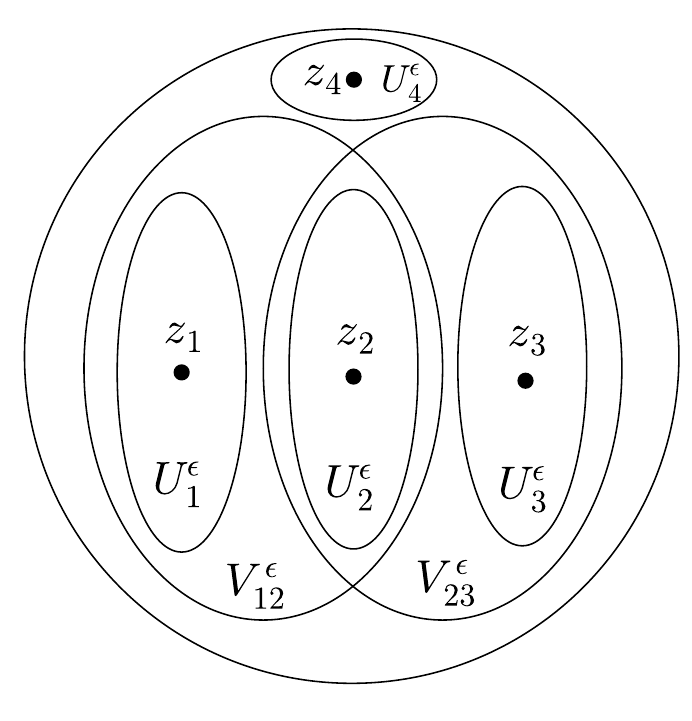}
  \caption{Open subsets in $\mathbb{S}$} \label{fig:sphere}
  \end{figure}
  We calculate $\gamma(\rho_{\epsilon}(\alpha),\x) \in K(\mathbb{S},4) $
  for $\alpha \in F_2$ and $\x \in X_4(\mathbb{S})$.
  We say that $\x =(x_1,x_2,x_3,x_4) \in X_4(\mathbb{S})$ is an \emph{$\epsilon$-good point} if all of $x_1$, $x_2$, $x_3$, and $x_4$ are in $U^{\epsilon}$.
We say that an $\epsilon$-good point $\bar{x}$ is \emph{of type $(p,q,r,s)$} if
 \begin{align*}
 \# (U^\epsilon_1 \cap \{ x_1, x_2, x_3, x_4\}) = p, \quad \# (U^\epsilon_2 \cap \{ x_1, x_2, x_3, x_4\}) = q, \\ \# (U^\epsilon_3 \cap \{ x_1, x_2, x_3, x_4\}) = r, 
\quad \# (U^\epsilon_4 \cap \{ x_1, x_2, x_3, x_4\}) = s.
 \end{align*}

Let $X^{\epsilon}_{pqrs}$ denote the set of $\epsilon$-good points $\x$ is of type $(p,q,r,s)$.
We define a cochain $c^{\epsilon}_{pqrs} \in C_b^n(F_2)$ by
\[  c_{pqrs}^{\epsilon}(\alpha_0,\dots,\alpha_n) = \int_{X^{\epsilon}_{pqrs}} c( \gamma(\rho_{\epsilon}(\alpha_0),\x),\dots, \gamma(\rho_{\epsilon}(\alpha_n),\x) )d\x\]
for $\alpha_0,\dots,\alpha_n \in F_2$.
By an argument similar to the proof of Lemma \ref{lem:key}, for $[c_{pqrs}^{\epsilon}]$ to be non-zero, both $W^{\epsilon}_{12}$ and $W^{\epsilon}_{23}$ must contain exactly two points, since the full twist of three or four strands is in the center $Z(P_4(\mathbb{S}))$. Thus, if $(p,q,r,s)$ is not $(1,1,1,1)$, $(0,2,0,2)$ or $(2,0,2,0)$, then $[c_{pqrs}^{\epsilon}]=0$.

Let $\x \in X_4(\mathbb{S})$ be an $\epsilon$-good point of type $(1,1,1,1)$, $(0,2,0,2)$ or $(2,0,2,0)$.
Similarly to Lemma \ref{lem:main_obs}, 
for almost every $\x \in X_4(\mathbb{S})$, there exists $\beta(\x) \in K(\mathbb{S},4)=Push(B_4(\mathbb{S}))$ such that
\[ \beta(\x)^{-1} \gamma(\rho_{\epsilon}(\alpha),\x) \beta(\x) =
  \begin{cases}
     \alpha & \text{type $(1,1,1,1)$}, \\
     (\overline{a_1})^{s_1(\alpha)+s_2(\alpha)} & \text{type $(0,2,0,2)$},\\
     (\overline{a_1})^{s_1(\alpha)}(\overline{a_3})^{s_2(\alpha)} & \text{type $(2,0,2,0)$},
  \end{cases}\]
for $\alpha \in F_2$.
Hence, we can prove $[c_{0202}^{\epsilon}]=[c_{2020}^{\epsilon}]=0$ and
\[ [c_{1111}^{\epsilon}] = \vol(X^{\epsilon}_{1111}) \cdot j^{\ast}(u) \]
by an argument similar to the proof of Lemma \ref{lem:key}.
Therefore,
\[\lim_{\epsilon \to 0} \| \rho_{\epsilon}^{\ast}(\REG_b^{\mathbb{S}}(u)) - \Lambda \cdot j^{\ast} (u) \| = 0\]
by setting
\[\Lambda = \lim_{\epsilon \to 0} 4!\cdot \area(U_1^{\epsilon})
\area(U_2^{\epsilon})\area(U_3^{\epsilon})\area(U_4^{\epsilon}). \qedhere\]
\end{proof}

\subsection{The case of a torus}
Let $\mathbb{T}$ denote the 2-torus. 
We summarize several facts on $B_2(\mathbb{T})$ (see \cite[Section 2.2]{BKS} for example). 
Recall that $\z=(z_1,z_2)$ denotes a base point of $X_2(\mathbb{T})$.
We define a braid $a_1$ (resp. $b_1$) so that it moves $z_1$ to the meridian (resp. longitude) direction and rotates once and does not move $z_2$.
We define $a_2$ and $b_2$ similarly by exchanging the role of $z_1$ and $z_2$.
It is known that $P_2(\mathbb{T}) \cong F_2 \times \ZZ^2$.
Namely, the set $\{\overline{\mathstrut a_1},\overline{\mathstrut b_1}\}$ generates $P_2(\mathbb{T})/Z(P_2(\mathbb{T})) \cong F_2$ and $\{a_1a_2,b_1b_2\}$ generates $Z(P_2(\mathbb{T})) \cong \ZZ^2$.
As in Section \ref{sec:sphere}, we obtain the map $j \colon F_2 \to K(\mathbb{T},2)$ induced by $Push \circ i \colon P_2(\mathbb{T}) \to K(\mathbb{T},2)$. Here, we consider
$P_2(\mathbb{T})/Z(P_2(\mathbb{T}))$ as $F_2=\langle \overline{\mathstrut a_1} , \overline{\mathstrut b_1} \rangle$.

\begin{lem} \label{lem:torus}
There exist a constant $\Lambda > 0$ and a family of homomorphisms $\{\rho_{\epsilon} \colon F_2 \to \G_{\mathbb{T}} \}_{0<\epsilon <1}$ such that
\[\lim_{\epsilon \to 0} \| \rho_{\epsilon}^{\ast}( \REG_b^{\mathbb{T}} (u) ) - \Lambda \cdot j^{\ast}(u)\| = 0 \]
for any $u \in \REH^n_b(K(\mathbb{T}),2)$.
\end{lem}

\begin{proof}
  For each $\epsilon$, we take open neighborhoods $U_i^{\epsilon}$ of $z_i \in \mathbb{T}$ $(i=1,2)$ so that
  \begin{itemize}
    \item $U_1^{\epsilon} \cap U_2^{\epsilon} = \emptyset$ and
    \item $\area(U^{\epsilon})=1-\epsilon$, where $U^{\epsilon}=U^{\epsilon}_1 \cup U^{\epsilon}_2$.
  \end{itemize}
  Moreover, we take open subsets $W^{\epsilon}_a$ and $V^{\epsilon}_a$ of $\mathbb{T}$ that are diffeomorphic to an annulus so that
  \begin{itemize}
    \item $U_1^{\epsilon} \subset \overline{W^{\epsilon}_{a}} \subset V^{\epsilon}_{a}$,
    \item $ V^{\epsilon}_a \cap  U_2^{\epsilon} = \emptyset$ and
    \item $W^{\epsilon}_{a}$ and $V^{\epsilon}_{a}$ contain a meridian.
  \end{itemize}
  We also take $W^{\epsilon}_{b}$ and $V^{\epsilon}_{b}$ similarly but to contain a longitude (see Figure \ref{fig:torus}).

\begin{figure}[t]
  \begin{minipage}[b]{0.4\linewidth}
      \centering
     \includegraphics[width=5cm]{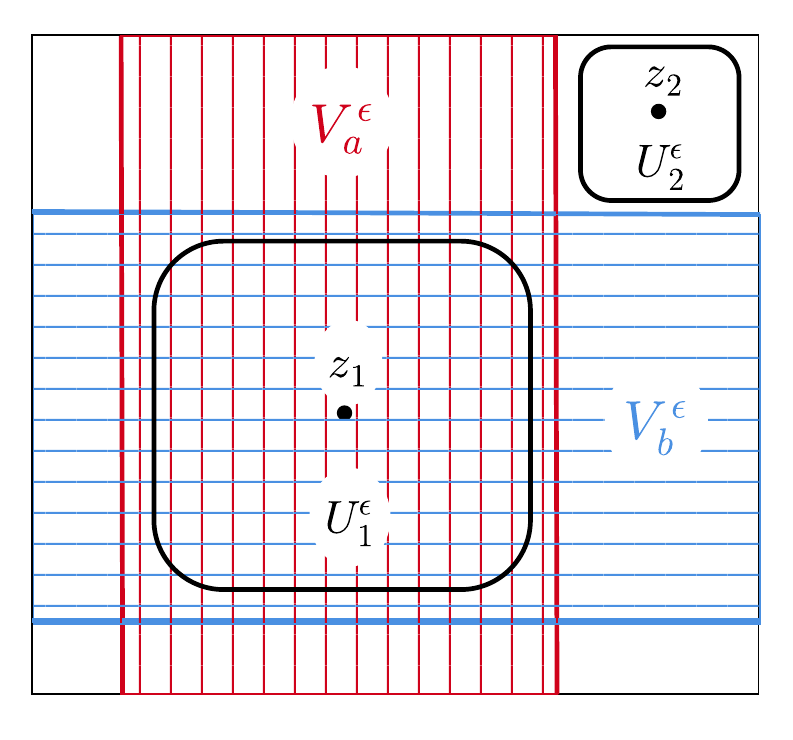}
   \end{minipage}
  \begin{minipage}[b]{0.5\linewidth}
    \centering
     \includegraphics[width=6cm]{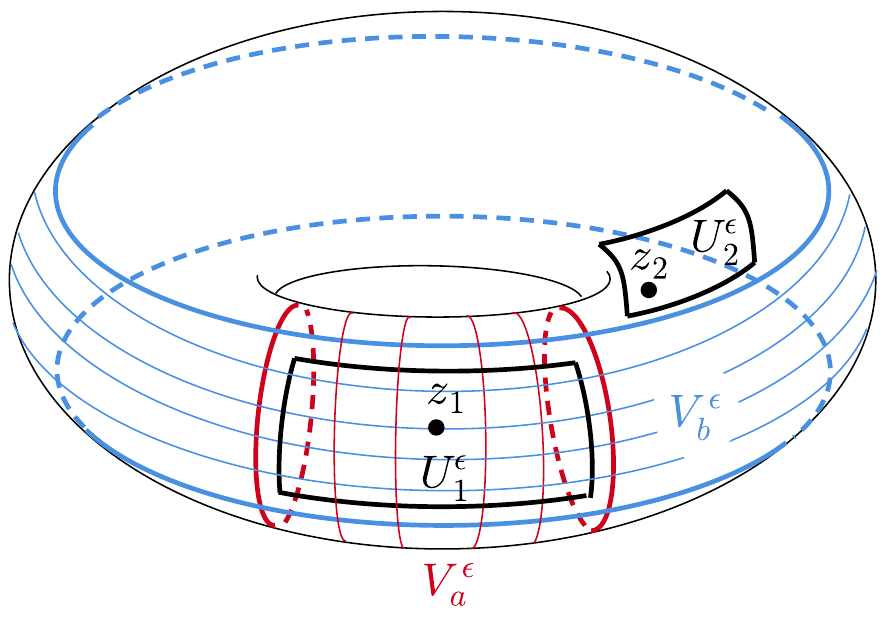}
   \end{minipage}
   \caption{Open subsets in $\mathbb{T}$}
    \label{fig:torus}
\end{figure}

We define $\rho_{\epsilon} \colon F_2\to \G_{\mathbb{T}}$ as follows.
We take an isotopy $\{(g_a)_t\}$ that rotates $W_a^{\epsilon}$ once and whose support is contained in $V_a^{\epsilon}$. For the generator $a_1 \in F_2$, we define $\rho_{\epsilon}(a_1)=(g_a)_1$. We also define $\rho_{\epsilon}(b_1)$ similarly.

We say that $\x =(x_1,x_2) \in X_2(\mathbb{T})$ is an \emph{$\epsilon$-good point} if both $x_1$ and $x_2$ are in $U^{\epsilon}$.
We say that an $\epsilon$-good point $\bar{x}$ is \emph{of type $(p,q)$} if
\[ \# (U^\epsilon_1 \cap \{ x_1, x_2\}) = p, \quad \# (U^\epsilon_2 \cap \{ x_1, x_2\}) = q. \]

Let $\x \in X_2(\mathbb{T})$ be an $\epsilon$-good point of type $(p,q)$.
We take an isotopy $\{(g_a)_t\}$ defined above.
Then, similar to Lemma \ref{lem:main_obs} again, 
for almost every $\x \in X_2(\mathbb{T})$, there exists $\beta'(\x) \in B_2(\mathbb{T})$ such that

\[
\gamma(\{(g_a)_t\},\x) =
\begin{cases}
  e & (p=0), \\
  \beta'(\x) a_1 \beta'(\x)^{-1} & (p=1),\\
  \beta'(\x) a_1a_2 \beta'(\x)^{-1} & (p=2).
\end{cases}
\]
Thus, we obtain  
\[
\gamma(\rho_{\epsilon}(\overline{a_1}),\x)= 
\begin{cases}
  \beta(\x) \overline{a_1}\beta(\x)^{-1} & (p=1),\\
  e & (\text{otherwise}), \\
\end{cases}
\]
where $\beta(\x)= Push(\beta'(\x)) \in K(\mathbb{T},2)$. Similarly, we can see that
\[
\gamma(\rho_{\epsilon}(\overline{b_1}),\x) =
\begin{cases}
  \beta(\x) \overline{b_1} \beta(\x)^{-1} & (q=1),\\
  e & (\text{otherwise}). \\
\end{cases}
\]
Hence, for $\alpha \in F_2$,
$\gamma(\rho_{\epsilon}(\alpha),\x)= \beta(\x) \alpha \beta(\x)^{-1} $ if $\x$ is of type
$(1,1)$.
By the argument similar to the proof of Lemma \ref{lem:key}, we can prove that
\[\lim_{\epsilon \to 0} \| \rho_{\epsilon}^{\ast}(\REG_b^{\mathbb{T}}(u)) - \Lambda \cdot  j^{\ast} (u) \| = 0\]
by setting
$\Lambda = \lim_{\epsilon \to 0} 2! \cdot \area(U^{\epsilon}_1) \area(U^{\epsilon}_2)$.
\end{proof}

\subsection{Remaining proofs}
We complete the proof of Theorems \ref{thm:S2T2}, \ref{thm:surface} and Corollary \ref{cor:S2T2}.

\begin{proof}[Proof of Theorem \ref{thm:S2T2}]
  By Lemmas \ref{lem:annulus}, \ref{lem:sphere}, and \ref{lem:torus}, we can prove Theorem \ref{thm:S2T2} by the same argument as in the proof of Theorem \ref{thm:main}.
\end{proof}

\begin{proof}[Proof of Corollary \ref{cor:S2T2}]

Since $B_2(\mathbb{A}) \cong K(\mathbb{A},2)$ is a finite index subgroup of $B_3$,
the inclusion map $K(\mathbb{A},2) \to B_3$ induces an isometric injective map $EH_b^3(B_3) \to EH_b^3(K(\mathbb{A},2))$ by Theorem \ref{thm:amenble_galois}. 
As we saw in the proof of Theorem \ref{thm:main}, $\REH_b^3(B_3)$ is uncountably infinite-dimensional, and thus 
$\REH_b^3( K(\mathbb{A},2))$ is also uncountably infinite-dimensional.

It is known that $\mathrm{MCG}(\mathbb{S},4)$ surjects onto $PSL(2,\ZZ)$ and its kernel is $\ZZ/2\ZZ \times \ZZ/2\ZZ$ (see \cite[Proposition 2.7]{MR2850125}).
Thus $\mathrm{MCG}(\mathbb{S},4)$ is quasi-isometric to $PSL(2,\ZZ)$.
Since $PSL(2,\ZZ)$ is non-elementary hyperbolic, $\mathrm{MCG}(\mathbb{S},4)$ is also.
Hence, by Theorems \ref{thm:acyl_hyp} and \ref{thm:kernel_forget}, $\REH_b^3( K(\mathbb{S},4) ) \cong \REH_b^3(\mathrm{MCG}(\mathbb{S},4) )$ is also uncountably infinite-dimensional.

Set $G = B_2(\mathbb{T})/Z(B_2(\mathbb{T}))$. Then $G$ has a presentation
\[ G = \langle a,b,c \mid a^2=b^2=c^2=1  \rangle \]
\cite[Exercise 6.3]{MuKu}.
Thus $G$ is isomorphic to $\ZZ/2\ZZ \ast \ZZ/2\ZZ \ast \ZZ/2\ZZ$ and hence it is non-elementary hyperbolic (because it is virtually free).
Hence, by Theorems \ref{thm:acyl_hyp} and \ref{thm:kernel_forget}, $\REH_b^3(G) \cong \REH_b^3( K(\mathbb{T},2) )$ is uncountably infinite-dimensional.

Therefore, we obtain the conclusion by Theorem \ref{thm:S2T2}.
\end{proof}

\begin{proof}[Proof of Theorem \ref{thm:surface}]
If $\chi(\Sigma) \geq 0$ then, by Corollary \ref{cor:S2T2}, $\REH_b^3(\G_{\Sigma})$ is infinite-dimensional. If $\chi(\Sigma) < 0$, $\pi_1(\Sigma)$ is a non-elementary hyperbolic group. Therefore, by the result of Brandenbursky and Marcinkowski \cite{brandenbursky2019bounded}, $\REH_b^3(\G_{\Sigma})$ is infinite-dimensional.
\end{proof}

\section*{Acknowledgment}
The author would like to thank Takahiro Matsushita for suggesting him to consider the case for the sphere and the torus after the author proved Theorem \ref{thm:main}.
He would like to thank Morimichi Kawasaki for informing him the paper
\cite{BKS}, which was helpful to consider the torus case. 
He is also grateful to both of them for reading the manuscript and providing comment to improve this paper.
He would like to thank Michael Brandenbursky, Micha{\l} Marcinkowski, and Shuhei Maruyama for helpful discussions.
The author is supported by JSPS KAKENHI Grant Number
JP20H00114 and JST-Mirai Program Grant Number JPMJMI22G1.

\bibliography{bddcoh_area_pres}

\providecommand{\bysame}{\leavevmode\hbox to3em{\hrulefill}\thinspace}
\providecommand{\MR}{\relax\ifhmode\unskip\space\fi MR }
\providecommand{\MRhref}[2]{%
  \href{http://www.ams.org/mathscinet-getitem?mr=#1}{#2}
}
\providecommand{\href}[2]{#2}
\begin{thebibliography}{10}

\bibitem{Birman}
J.~S. Birman, \emph{Braids, links, and mapping class groups}, Princeton
  University Press, Princeton, N.J.; University of Tokyo Press, Tokyo, 1974,
  Annals of Mathematics Studies, No. 82.

\bibitem{Bra_knot}
M.~Brandenbursky, \emph{On quasi-morphisms from knot and braid invariants}, J.
  Knot Theory Ramifications \textbf{20} (2011), no.~10, 1397--1417.

\bibitem{BKS}
M.~Brandenbursky, J.~K{\k{e}}dra, and E.~Shelukhin, \emph{On the autonomous
  norm on the group of {H}amiltonian diffeomorphisms of the torus}, Commun.
  Contemp. Math. \textbf{20} (2018), no.~2, 1750042, 27.

\bibitem{BM_ent}
M.~Brandenbursky and M.~Marcinkowski, \emph{Entropy and quasimorphisms}, J.
  Mod. Dyn. \textbf{15} (2019), 143--163.

\bibitem{brandenbursky2019bounded}
\bysame, \emph{Bounded cohomology of transformation groups}, Math. Ann.
  \textbf{382} (2022), no.~3-4, 1181--1197.

\bibitem{MR2850125}
B.~Farb and D.~Margalit, \emph{A primer on mapping class groups}, Princeton
  Mathematical Series, vol.~49, Princeton University Press, Princeton, NJ,
  2012.

\bibitem{FPS}
R.~Frigerio, M.~B. Pozzetti, and A.~Sisto, \emph{Extending higher-dimensional
  quasi-cocycles}, J. Topol. \textbf{8} (2015), no.~4, 1123--1155.

\bibitem{GG}
J.-M. Gambaudo and \'{E}. Ghys, \emph{Commutators and diffeomorphisms of
  surfaces}, Ergodic Theory Dynam. Systems \textbf{24} (2004), no.~5,
  1591--1617.

\bibitem{GP}
J.-M. Gambaudo and E.~E. P{\'{e}}cou, \emph{Dynamical cocycles with values in
  the {A}rtin braid group}, Ergodic Theory Dynam. Systems \textbf{19} (1999),
  no.~3, 627--641.

\bibitem{Gromov}
M.~Gromov, \emph{Volume and bounded cohomology}, Inst. Hautes \'{E}tudes Sci.
  Publ. Math. (1982), no.~56, 5--99 (1983).

\bibitem{GJP}
J.~Guaschi and D.~Juan-Pineda, \emph{A survey of surface braid groups and the
  lower algebraic {$K$}-theory of their group rings}, Handbook of group
  actions. {V}ol. {II}, Adv. Lect. Math. (ALM), vol.~32, Int. Press,
  Somerville, MA, 2015, pp.~23--75.

\bibitem{HD}
M.-E. Hamstrom and E.~Dyer, \emph{Regular mappings and the space of
  homeomorphisms on a 2-manifold}, Duke Math. J. \textbf{25} (1958), 521--531.

\bibitem{Ish}
T.~Ishida, \emph{Quasi-morphisms on the group of area-preserving
  diffeomorphisms of the 2-disk via braid groups}, Proc. Amer. Math. Soc. Ser.
  B \textbf{1} (2014), 43--51.

\bibitem{KT}
C.~Kassel and V.~Turaev, \emph{Braid groups}, Graduate Texts in Mathematics,
  vol. 247, Springer, New York, 2008, With the graphical assistance of Olivier
  Dodane.

\bibitem{MR1902362}
R.~P. Kent, IV and D.~Peifer, \emph{A geometric and algebraic description of
  annular braid groups}, vol.~12, 2002, International Conference on Geometric
  and Combinatorial Methods in Group Theory and Semigroup Theory (Lincoln, NE,
  2000), pp.~85--97.

\bibitem{MuKu}
K.~Murasugi and B.~I. Kurpita, \emph{A study of braids}, Mathematics and its
  Applications, vol. 484, Kluwer Academic Publishers, Dordrecht, 1999.

\bibitem{Nit}
M.~Nitsche, \emph{Higher-degree bounded cohomology of transformation groups},
  preprint, arXiv:2105.08698 (2021).

\bibitem{Osin}
D.~Osin, \emph{Acylindrically hyperbolic groups}, Trans. Amer. Math. Soc.
  \textbf{368} (2016), no.~2, 851--888.

\bibitem{Smale}
S.~Smale, \emph{Diffeomorphisms of the {$2$}-sphere}, Proc. Amer. Math. Soc.
  \textbf{10} (1959), 621--626.

\bibitem{Va}
V.~A. Vassiliev, \emph{Complements of discriminants of smooth maps: topology
  and applications}, Translations of Mathematical Monographs, vol.~98, American
  Mathematical Society, Providence, RI, 1992, Translated from the Russian by B.
  Goldfarb.

\end{thebibliography}
\bibliographystyle{amsplain}

\end{document}